\documentclass[12pt]{amsart}
\usepackage{amsfonts}

\theoremstyle{definition}

\theoremstyle{remark}

\newcommand{\const}{\mathop{\rm const}\limits}

\newcommand{\meas}{\mathop{\rm meas}\limits}

\newcommand{\grad}{\mathop{\rm grad}\limits}

\begin{document}

\begin{center}

{\bf POINCAR\'E INEQUALITIES IN BILATERAL GRAND LEBESGUE SPACES} \par

\vspace{3mm}
{\bf E. Ostrovsky}\\

e - mail: galo@list.ru \\

\vspace{3mm}

{\bf L. Sirota}\\

e - mail: sirota@zahav.net.il \\

\vspace{3mm}

{\bf E.Rogover}\\

e - mail: rogovee@gmail.com

\vspace{4mm}

 Abstract. \\
{\it In this paper we obtain the non - asymptotic estimations of Poincare type
between function  and its gradient in the so - called
Bilateral Grand Lebesgue Spaces. We also give some examples to show the
sharpness of these inequalities.} \\

\end{center}

\vspace{3mm}


\vspace{3mm}


\vspace{3mm}

\section{Introduction}

\vspace{3mm}

"The term Poincar\'e type inequality is used, somewhat loosely,
to describe a class of inequalities that generalize the classical Poincar\'e
inequality"

$$
\int_D |f(x)|^p \ dx \le A(p,D) \int_D | \ |\grad f(x)|^p \ | \ dx, \eqno(0)
$$
see \cite{Adams1}, chapter 8,  p.215, and the source work of Poincar\'e
\cite{Poincare1}. \par
 We will call "the Poincar\'e inequality" some {\it improved, or modified }
 inequality

$$
 |T_{\delta} f(\cdot)|_p  \le C(D) \cdot \frac{p}{|p - d|} \cdot
 | \ |\grad f(\cdot)| \ |_p, \eqno(1)
$$
or  more generally the inequality of a view:

$$
 |T_{\delta} f^0(\cdot)/(\delta(x))^{\alpha}|_p  \le C_{\alpha}(D) \cdot
  \frac{p}{|p - d(1 + \alpha)|} \cdot | \ |\grad f(x) | /\delta(x)^{\alpha} \ |_p. \eqno(2)
$$
and we will call the domains $ \{ D \}, \ D \subset R^d, \ d = 2,3,\ldots $
which satisfied the inequalities (1) or (2) for any functions
$ f $ in the Sobolev class $ W(1,p) $  as a {\it Poincar\'e domains.} \par
 There  are many publications about these inequalities and its applications, see,
 for instance, \cite{Chuas1}, \cite{Fazio1}, \cite{Kufner1}, \cite{Rjtva1},
\cite{Wannebo1} and the classical monographs \cite{Beckenbach1},  \cite{Hardy1};
see also reference therein.\par
 Here
 $$
 \delta = \delta(x) = \delta_D(x) = \inf_{y \in \partial D} |x - y|
 $$
is the distance between the point $ x, \ x \in D $ and the boundary $ \partial D $
of the set $ D, \ p \in [1, \infty), \ p \ne d(1 + \alpha) $ in the case when the
domain $ D $ is bounded and $ \delta(x) = |x| $ otherwise;

$$
\alpha = \const \in (-1, \infty); \  x \in R^d \ \Rightarrow |x| = (x,x)^{1/2};
$$
following, for the real value $ x  \ |x| \ $ denotes usually absolute value of $ x;$ \par
 $$
 T_{\delta}f(x) = \frac{f^0(x)}{\delta(x)}, \ T_{\alpha,\delta}f(x) =
 \frac{f^0(x)}{(\delta(x))^{1 + \alpha}} =
 T_{\delta}f(x)/\left[(\delta(x))^{\alpha}\right],
 $$
where
$$
f^0(x) = f(x) - \int_D f(x) dx /|D|, \ |D| = \meas(D) \in (0,\infty)
$$
in the case of bounded domain $ D $ (a "centering" of a function $ f,) $
and $ f^0(x) = f(x) $ otherwise.\par
 Note that in the case of bounded domain $ D $ the inequalities (1) and (2) may be
 rewritten as follows:

 $$
 \inf_{c \in R} | (f(x) - c)/(\delta(x))^{1 + \alpha} |_p  \le C_{\alpha}(D) \cdot
\frac{p}{|p - d(1 + \alpha)|} \cdot | \ |\grad f(x) | /\delta(x)^{\alpha} \ |_p.
 $$

 We will called the operators $ T_{\delta} $ and $ T_{\alpha,\delta} $  as a
Poincar\'e operators. \par
 It is known that if the domain $ D $ is open, bounded, contain the origin and has
 a Lipschitz or at last H\"older  boundary, or consists on the finite union of these domains, that it is Poincar\'e domain. \par

  \vspace{3mm}

{\sc We will assume that the considered set $ D $ is Poincar\'e domain, at last
in some parameter } $ \alpha. $ \par

\vspace{3mm}

{\sc We will distinguish a two cases: the first, or bounded case is when the domain
$ D $ is bounded and satisfied the Poincar\'e condition, and the second case of
unbounded domain also with Poincar\'e property.}\par

\vspace{3mm}

In the first case the value $ p $ might belongs to the semi-closed interval
$  1 \le p < d $ or in more general case when $ \alpha \ne 0 $ we suppose
 $ d(1 + \alpha) > 1 $ and in the second case or correspondingly
 $ d(1 + \alpha) < p < \infty. $  \par
Note that in the case of the bounded domain $ D $ the Poincar\'e inequality (1) may
be rewritten as follows:

$$
 |T_{\delta} f^0(x)/\delta(x)^{\alpha}|_p  \le C_{\alpha}(D) \cdot
  \frac{1}{|p - d(1 + \alpha)|} \cdot | \ |\grad f(x) | /\delta(x)^{\alpha} \ |_p, \eqno(3)
$$

 We denote as usually the classical $ L_p $ Lebesgue norm

 $$
 |f|_p = |f|_{p,D} = \left( \int_{D} |f(x)|^p \ dx  \right)^{1/p};  \  f \in L_p \
 \Leftrightarrow |f|_p < \infty,
 $$
and denote for arbitrary measurable subset $ D_1 $ of the set $ D: \ D_1 \subset D $

$$
|f|_{p,D_1} = \left( \int_{D_1} |f(x)|^p \ dx  \right)^{1/p};
$$
and  denote also $ L(a,b) = \cap_{p \in (a,b)} L_p.  $ \par

{\bf Lemma 1.} If $ D_1, D_2 \subset D, \ D_1 \cap D_2 = \emptyset, $ then

$$
|f|_{p,D_1 \cup D_2} \le |f|_{p,D_1} + |f|_{p,D_2}.
$$
The proof follows immediately from the triangle inequality for the $ L_p $ norm.\par
{\bf Lemma 2.} Let $ D $ be a measurable set with
non-trivial finite Lebesgue  measure: $ 0 < \mu(D) < \infty. $ Let also $ f(x),
f: D \to R $ be positive a.e. bounded function. We assert that for all the values
$ p \in [1,\infty) $

$$
0 < C_1 \le |f|_p \le C_2 < \infty.
$$
{\bf Proof.} \\
{\bf 1. Upper bound.} Let us denote $ M = \sup_{x \in D} f(x). $ We have:

$$
|f|^p_p \le M^p \cdot \mu(D),
$$
following

$$
|f|_p \le M \cdot [\mu(D)]^{1/p} \le M \ \max(1, \mu(D)) \stackrel{def}{=} C_2.
$$

{\bf 2. Low bound.}  Let $ D_1 $ be some subset of the set $ D $ with positive measure
and $ m $ be some positive number such that

$$
\forall p \in D_1 \ f(x) \ge m.
$$
We have:
$$
|f|^p_p \ge m^p \cdot  \mu(D_1),
$$
following
$$
|f|_p \ge m \cdot \left[ \mu(D_1) \right]^{1/p} \ge m \cdot \min(1, \mu(D_1))
\stackrel{def}{=} C_1.
$$

\vspace{3mm}

{\it  Our aim is a generalization of the estimation (2), (3) on the so - called Bilateral
 Grand Lebesgue Spaces $ BGL = BGL(\psi) = G(\psi), $ i.e. when } $ f(\cdot)
 \in G(\psi) \ $  {\it and to show the precision of obtained estimations by means
 of the constructions of suitable examples. } \par

 \vspace{3mm}

  We recall briefly the definition and needed properties of these spaces.
  More details see in the works \cite{Fiorenza1}, \cite{Fiorenza2}, \cite{Ivaniec1},
   \cite{Ivaniec2}, \cite{Ostrovsky1}, \cite{Ostrovsky2}, \cite{Kozatchenko1},
  \cite{Jawerth1}, \cite{Karadzov1} etc. More about rearrangement invariant spaces
  see in the monographs \cite{Bennet1}, \cite{Krein1}. \par

\vspace{3mm}

For $a$ and $b$ constants, $1 \le a < b \le \infty,$ let $\psi =
\psi(p),$ $p \in (a,b),$ be a continuous positive
function such that there exists a limits (finite or not)
$ \psi(a + 0)$ and $\psi(b-0),$  with conditions $ \inf_{p \in (a,b)} > 0 $ and
 $\min\{\psi(a+0), \psi(b-0)\}> 0.$  We will denote the set of all these functions
 as $ \Psi(a,b). $ \par

The Bilateral Grand Lebesgue Space (in notation BGLS) $  G(\psi; a,b) =
 G(\psi) $ is the space of all measurable
functions $ \ f: R^d \to R \ $ endowed with the norm

$$
||f||G(\psi) \stackrel{def}{=}\sup_{p \in (a,b)}
\left[ \frac{ |f|_p}{\psi(p)} \right], \eqno(4)
$$
if it is finite.\par
 In the article \cite{Ostrovsky2} there are many examples of these spaces.
 For instance, in the case when  $ 1 \le a < b < \infty, \beta, \gamma \ge 0 $ and

 $$
 \psi(p) = \psi(a,b; \beta, \gamma; p) = (p - a)^{-\beta} (b - p)^{-\gamma};
 $$
we will denote
the correspondent $ G(\psi) $ space by  $ G(a,b; \beta, \gamma);  $ it
is not trivial, non - reflexive, non - separable
etc.  In the case $ b = \infty $ we need to take $ \gamma < 0 $ and define

$$
\psi(p) = \psi(a,b; \beta, \gamma; p) = (p - a)^{-\beta}, p \in (a, h);
$$

$$
\psi(p) = \psi(a,b; \beta, \gamma; p) = p^{- \gamma} = p^{- |\gamma|}, \ p \ge h,
$$
where the value $ h $ is the unique  solution of a continuity equation

$$
(h - a)^{- \beta} = h^{ - \gamma }
$$
in the set  $ h \in (a, \infty). $ \par

 The  $ G(\psi) $ spaces over some measurable space $ (X, F, \mu) $
with condition $ \mu(X) = 1 $  (probabilistic case)
appeared in \cite{Kozatchenko1}.\par
 The BGLS spaces
are rearrangement invariant spaces and moreover interpolation spaces
between the spaces $ L_1(R^d) $ and $ L_{\infty}(R^d) $ under real interpolation
method \cite{Carro1}, \cite{Jawerth1}. \par
It was proved also that in this case each $ G(\psi) $ space coincides
with the so - called {\it exponential Orlicz space,} up to norm equivalence. In others
quoted publications were investigated, for instance,
 their associate spaces, fundamental functions
$\phi(G(\psi; a,b);\delta),$ Fourier and singular operators,
conditions for convergence and compactness, reflexivity and
separability, martingales in these spaces, etc.\par

{\bf Remark 1.} If we introduce the {\it discontinuous} function

$$
\psi_r(p) = 1, \ p = r; \psi_r(p) = \infty, \ p \ne r, \ p,r \in (a,b)
$$
and define formally  $ C/\infty = 0, \ C = \const \in R^1, $ then  the norm
in the space $ G(\psi_r) $ coincides with the $ L_r $ norm:

$$
||f||G(\psi_r) = |f|_r.
$$

Thus, the Bilateral Grand Lebesgue spaces are direct generalization of the
classical exponential Orlicz's spaces and Lebesgue spaces $ L_r. $ \par

The BGLS norm estimates, in particular, Orlicz norm estimates for
measurable functions, e.g., for random variables are used in PDE
\cite{Fiorenza1}, \cite{Ivaniec1}, theory of probability in Banach spaces
\cite{Ledoux1}, \cite{Kozatchenko1},
\cite{Ostrovsky1}, in the modern non-parametrical statistics, for
example, in the so-called regression problem \cite{Ostrovsky1}.\par

The article is organized as follows. In the next section we obtain
the main result: upper bounds for Poincar\'e  operators in the Bilateral
Grand Lebesgue spaces. In the third section we construct some examples
in order to illustrate the precision of upper estimations. \par

 The last section contains some slight generalizations of obtained results.\par

\vspace{3mm}

 We use symbols $C(X,Y),$ $C(p,q;\psi),$ etc., to denote positive
constants along with parameters they depend on, or at least
dependence on which is essential in our study. To distinguish
between two different constants depending on the same parameters
we will additionally enumerate them, like $C_1(X,Y)$ and
$C_2(X,Y).$ The relation $ g(\cdot) \asymp h(\cdot), \ p \in (A,B), $
where $ g = g(p), \ h = h(p), \ g,h: (A,B) \to R_+, $
denotes as usually

$$
0< \inf_{p\in (A,B)} h(p)/g(p) \le \sup_{p \in(A,B)}h(p)/g(p)<\infty.
$$
The symbol $ \sim $ will denote usual equivalence in the limit
sense.\par
We will denote as ordinary the indicator function
$$
I(x \in A) = 1, x \in A, \ I(x \in A) = 0, x \notin A;
$$
here $ A $ is a measurable set.\par
 All the passing to the limit in this article may be grounded by means
 of Lebesgue dominated convergence theorem.\par

\bigskip

\section{Main result: upper estimations for Poincar\'e operator}

\vspace{3mm}

Let $ \psi(\cdot) \in \Psi(a,b), \ $ where $ 1 = a < b = d $ in the case of bounded
domain $ D $ and $ d = a < b = \infty $ otherwise. \par
 Define for the arbitrary function  $ \psi \in \Psi(a,b)  $ the auxiliary
 function of the {\it variable} $ \ p $

$$
\psi_{\alpha, d}(p) =\frac{p}{| d - p(1 + \alpha)|} \cdot \psi(p), \eqno(4)
$$

$$
\psi_{d}(p) =\frac{p}{| d - p|} \cdot \psi(p) = \psi_{0,d}(p).
$$

Notice that only the values $ p = \infty $ and $ p = d(1 + \alpha) $ are
critical points; another points are not interest. \par

\vspace{3mm}

{\bf Theorem 1.} Let $ f \in G(\psi), \ \psi \in \Psi(a,b) $  and let the domain
$ D $ be the Poincar\'e domain.  Then

\vspace{3mm}

$$
||T_{\alpha, \delta} \ f||G(\psi_{\alpha, d}) \le C( D,\alpha, d) \
|| \ | \grad f \ | \ ||G(\psi).  \eqno(5).
$$

\vspace{3mm}

{\bf Example 1.} When the domain $ D $ is bounded (the first case),
 $ a = 1, b = d,  \beta, \gamma > 0, $ and
$  f \in G(1,d; \beta, \gamma), \ f \ne 0, $ then

$$
T_{\delta} f(\cdot) \in G(1, d; \beta, \gamma + 1 ). \eqno(6)
$$

\vspace{3mm}

{\bf Example 2.} When the domain $ D $ is  Poincar\'e and
unbounded (the second case),
 $ a = d, b = \infty,  \beta >0, \gamma < 0, $ and
$  f \in G(d, \infty; \beta, \gamma), \ f \ne 0, $ then

$$
T_{\delta} f(\cdot) \in G( d, \infty; \beta + 1, \gamma + 1 ). \eqno(7)
$$

\vspace{3mm}

{\bf Proof} of the theorem 1 is very simple. Denote for the simplicity
$ u = T_{\alpha, \delta}f^0; \ u: D \to R. $\par
 We suppose $ f(\cdot) \in G(\psi); $ otherwise is nothing to prove.\par

 We can assume without loss of generality that $ ||\  |\grad f| \ ||G(\psi) = 1; $
this means that

$$
\forall p \in (a,b) \ \Rightarrow |\ | \grad f| \ |_p \le \psi(p).
$$

Using the inequality (1) we obtain:

$$
|u|_p \le C(D) \frac{p}{|d - p(1 + \alpha)|} \cdot  \psi(p)  = C(D) \
\psi_{\alpha,d}(p) = C(D) \ \psi_{\alpha,d}(p) \ || \ |\grad f| \ ||G(\psi).
$$

The assertion of theorem 1 follows after
the dividing over the  $ \psi_{\alpha,d}(p), $ tacking the supremum over
$ p, \ p \in (a,b) $ and on the basis of the definition of the $ G(\psi) $ spaces.
\hfill $\Box$ \\

\bigskip

\section{Low bounds for Poincar\'e inequality.}

\vspace{3mm}

In this section we built some examples in order to illustrate the
 exactness of upper estimations. We consider both the cases
$ a = 1, \ b = d $  and $ a = d, \ b = \infty. $ \par

\vspace{3mm}

{\bf A. Bounded domain.}\\

\vspace{3mm}
Note that in the bounded case only the value $ p = p_0 \stackrel{def}{=}
d(1 + \alpha) $  is the critical value.  It is presumed in this (sub)section
that $ d(1 + \alpha) > 1; $ otherwise is nothing to prove. \par
 Let us denote for the mentioned  values $ p \in (1,d(1 + \alpha)) $ and for the
 function  $ \ f \in \cap_{p \in (1,d(1 + \alpha))} L_p, f \ne 0 $ the quantity
$$
V_{\alpha,d}(f,p) = V(f,p) = \frac{|T_{\alpha,d} f |_p \cdot
[(d - p(1 + \alpha))/p]}{|f|_p }, \eqno(8)
$$

From the inequality (2) follows that for {\it some } non - trivial (positive and
finite) number $ C^{(1)} = C^{(1)}(\alpha,d) $

$$
 \sup_{f \in L(1,d(1+\alpha)), f \ne 0} \overline{\lim}_{p \to d(1 +\alpha) - 0}
 V( f,p) \le C^{(1)}.
$$

We intend to prove an inverse inequality at the critical point
$ p \to d(1 + \alpha) - 0. $  \par

\vspace{3mm}

{\bf Theorem 2.a.}  For all the values $ \alpha \in (-1, \infty), $ such that
$ d(1 + \alpha) > 1 \ $ there exists a constants
$ C_1 = C_1(\alpha,d) \in (0,\infty) $ and the bounded Poincar\'e domain $ D $
for which

\vspace{3mm}

$$
 \sup_{f \in L(1,d(1+\alpha)), f \ne 0} \underline{\lim}_{p \to d(1+\alpha)-0}
 V_{\alpha,d}(f,p) \ge C_1. \eqno(9)
$$

\vspace{3mm}

{\bf Proof.}
  Let us consider a function (more exactly, a family of the functions) of a view

$$
u(x) = u_{\Delta}(x) = | \ \log |x| \ |^{\Delta} \cdot I(|x| < 1/e), \
\Delta = \const > 1.
$$
Here the domain $ D $ is the unit ball of the space $ R^d: \ D = \{x: |x| \le 1   \}. $
\par
Note that the average value $ |D|^{-1} \int_D u(x) \ dx $ is bounded.\par
We have using the multidimensional spherical coordinates as
$ p \to  d(1 + \alpha) - 0:$

$$
|(\delta(x))^{-\alpha} \ \grad u|^p_p \asymp
\int_0^{1/e} r^{d - 1 - p(1 + \alpha)} |\log r|^{p(\Delta - 1)} \ dr
$$

$$
\asymp \int_0^{\infty} \exp(-y(d - p(1+\alpha))) \ y^{p(\Delta - 1)} \ dy
$$

$$
= [d - p(1+\alpha)]^{-p(\Delta - 1) - 1} \ \Gamma(p(\Delta -1))
$$

$$
\asymp [d - p(1+\alpha)]^{-p(\Delta - 1) - 1}, \eqno(10)
$$
since we can suppose that the value $ p $ is bounded. Here the $ \Gamma(\cdot) $
denotes usually Gamma function. \par

 It follows from the equality (10) that

 $$
|(\delta(x))^{-\alpha} \ \grad u|_p \asymp
 [d - p(1 + \alpha) ]^{-(\Delta - 1) -(1 + \alpha)/d }.\eqno(11)
 $$

Further, we find analogously:

$$
| \delta^{-1-\alpha} u^0 |_p^p \asymp \int_1^{\infty} e^{-y(d - p(1+\alpha))} \
y^{p \Delta} \ dy \asymp [d - p(1+\alpha)]^{-p \Delta - 1},
$$
therefore

$$
| \delta^{-1-\alpha} u^0 |_p  \asymp [d - p(1 + \alpha)]^{- \Delta -(1 + \alpha)/d }.
\eqno(12)
$$

We conclude substituting the estimations  (11) and (12) into the expression for
the value $ V(f,p) $ that it is bounded from below as $ p \to d(1+\alpha) - 0.$ \par

But the function $ u = u_{\Delta}(x) $ is discontinuous. Let us redefine the function
$ u = u(x) $ as following. Let $ \phi = \phi(r) $ be infinitely differentiable
 function with support on the set $ r \in [1/e, 2/e]: $

 $$
 \phi(r) > 0 \ \Leftrightarrow r \in (1/e, 2/e)
 $$
and such that

$$
\phi(1/e) = 1, \  \phi^/(1/e) = e \cdot \Delta, \ \phi(2/e) = \phi^/(2/e) = 0.
$$
Then the function
$$
\tilde{u}_{\Delta}(x) = |\log |x| \ |^{\Delta} \cdot I(|x| \in (1, 1/e) + \phi(|x|) \cdot
I(|x| \in [1/e, 2/e])
$$

$$
\stackrel{def}{=} u_{\Delta}(x) +  w(x) = u_{\Delta}(x) \ I(|x| \in (0,1/e) +
w(x) \ I(|x| \in (1/e, 2/e)).
$$
gives us the example for the theorem 2.a. Indeed, since the supports for the functions
$ u_{\Delta} $ and $ w(x) $ are adjoint, we have using the lemma 1
for the nominator for the expression of a function $ V \ $ (8) the following simple estimation {\it from below:}

$$
\left| \tilde{u}_{\Delta} \right|_p  \ge | u_{\Delta} |_p.
$$
We need further to obtain the {\it upper estimate } for
the denominator for the expression for the function $ V. $  We have using the assertions
of lemma 1 and lemma 2:

$$
\left| \grad \tilde{u}_{\Delta}  \right|_p  \le \left| \grad u_{\Delta} \right|_p +
| \grad w |_p
$$

 $$
 \le C_5 \left[ \frac{p}{|p(1+\alpha) - d|}\right]^{(\Delta - 1) -(1 + \alpha)/d } +
 C_6
 $$

 $$
\le  C_7 \left[ \frac{p}{|p(1+\alpha) - d|}\right]^{(\Delta - 1) -(1 + \alpha)/d }.
 $$
This completes the proof of theorem 2.a. \par

\vspace{3mm}

{\bf B. Unbounded domain.}\\

\vspace{3mm}

{\bf Theorem 2.b.}  For all the values $ \alpha \in (-1, \infty), $ such that
$ d(1 + \alpha) > 1 \ $ there exists a constants
$ C_2 = C_2(\alpha,d) \in (0,\infty) $ and the unbounded Poincar\'e domain $ D $
for which

\vspace{3mm}

$$
 \sup_{f \in L(1,d(1+\alpha)), f \ne 0} \underline{\lim}_{p \to d(1+\alpha)+0}
 V_{\alpha,d}(f,p) \ge C_2. \eqno(13)
$$
and

$$
 \sup_{f \in L(1,d(1+\alpha)), f \ne 0} \underline{\lim}_{p \to \infty }
 V_{\alpha,d}(f,p) \ge C_2. \eqno(14)
$$

{\bf Proof} is at the same as in the proof of theorem 2.a. It is sufficient to
consider the domain $ D = \{x: \ |x| \ge 1 \} $ and to choose the {\it single}
function

$$
v(x) = v_{\Delta}(x) = | \ (\log |x|) \ |^{\Delta} \ I(|x| > e)+ I(0.5e \le |x| \le e )
\phi(|x|),
$$
where $ \phi = \phi(r) $ is infinitely differentiable function with the support
$ r \in (0.5 e, e) $ such that

$$
\phi(e/2) = \phi^/(e/2) = 0, \phi(e) = 1, \phi^/(e) = \Delta/e;
$$
$ \ \Delta = \const > 1.$ \par

We have for the values $ p \in ( d(1+\alpha) + 0, \infty) $ repeating the considerations
for the proof of theorem 2.a used the lemmas  1.2:

$$
| \ | \grad v| \ |x|^{-\alpha} \ |^p_p \asymp C_7^p \int_e^{\infty} r^{d-1-p(1+\alpha)} \
(\log r)^{p(\Delta - 1)} \ dr + C_8^p
$$

$$
\asymp C_7^p \int_0^{\infty} \exp(-y(p(1+\alpha) - d) \ y^{p(\Delta-1)} \ dy + C_8^p
= C_7^p \frac{\Gamma(p(\Delta - 1) + 1)}{[p(1 + \alpha) - d]^{p(\Delta-1) + 1}}
+ C_8^p;
$$

$$
| \ | \grad v_{\Delta}| \ |x|^{-\alpha} \ |_p  \asymp C_9 \ 
\frac{p^{\Delta - 1}}{[d - p(1+\alpha)]^{\Delta - 1 + 1/p }} + C_{10}
$$

$$
\asymp C_{11} \ \frac{p^{\Delta - 1}}{[d - p(1+\alpha)]^{\Delta - 1 + 1/p }};
$$

$$
\left| v_{\Delta}/|x|^{1 + \alpha} \right|^p_p \asymp C_{12}^p \ 
\frac{\Gamma(p\Delta + 1)}{[p(1 + \alpha) - d]^{p \Delta +1} };
$$

$$
\left| v_{\Delta}/|x|^{1 + \alpha} \right|_p \asymp C_{13} \ 
\frac{p^{\Delta }}{[d - p(1+\alpha)]^{\Delta + 1/p }}
$$

$$
\asymp C_{14} \ \frac{p^{\Delta }}{[d - p(1+\alpha)]^{\Delta + 1/p }};
$$
therefore

$$
\left[\left| v_{\Delta}/|x|^{1 + \alpha} \right|_p \right]:
\left[| \ | \grad v_{\Delta}| \ |x|^{-\alpha} \ |_p  \right] \asymp
p/[ d - p(1 + \alpha)].
$$

This completes the proof of theorem 2.b. We omitted  the simple calculations
used the Stirling's formula as $ p \to \infty $ etc. \par

\vspace{3mm}

 {\bf Remark 2.} It follows from theorems 2.a and 2.b that the estimation of theorem 1
 is exact in the $ G(\psi) $ spaces. Namely, for the function $ \psi_{\Delta}(p) $
 of a view

 $$
 \psi_{\Delta}(p) = |u_{\Delta}(\cdot) |_p
 $$
or correspondingly

 $$
 \psi_{\Delta}(p) = |v_{\Delta}(\cdot) |_p
 $$
in the assertion of theorem 1 (5) take place the inverse asymptotical inequality,
up to multiplicative constant.\par \hfill $\Box$

\bigskip

\section{Concluding remarks}

\vspace{3mm}
{\it We consider in this section some generalizations of inequalities (1),(2).} \par
{\bf A.} Let us denote
$$
\log^+ z = \max(1, |\log z| ), z > 0; \ z_+ = \max(z,0); \ z \in (-\infty, \infty).
$$
We obtain using the works
\cite{Chuas1}, \cite{Fazio1}, \cite{Kufner1}, \cite{Rjtva1}, \cite{Wannebo1}
etc. the following inequality

$$
\left| \frac{u^0 \cdot (\log^+ \delta)^{B_1}}{\delta^{1 + \alpha}} \right|_p \le
C(d,\alpha,D,B_1,B_2) \times
$$

$$
\left[ \frac{p}{d - p(1 + \alpha)} \right]^{(1 - B_2 + B_1)_+ } \cdot
\left| \frac{|\grad u| \cdot (\log^+ \delta)^{B_2}}{\delta^{\alpha}} \right|_p.
\eqno(15)
$$
Here in addition $ B_1, B_2 = \const, \ \min(B_1,B_2) > 0. $ \par
We assert analogously to the theorems 2.a,2.b and using at the same examples
that the inequality (15) is asymptotically exact as $ p \to d(1 + \alpha) $ and as
$ p \to \infty. $ \par

At the same result is true for more general operator:

$$
\left| \frac{u^0 \cdot (\log^+ \delta)^{B_1} \cdot S(\log^+ \delta )}
{\delta^{1 + \alpha}} \right|_p \le
C(d,\alpha,D,B_1,B_2) \times
$$

$$
\left[ \frac{p}{d - p(1 + \alpha)} \right]^{(1 - B_2 + B_1)_+ } \cdot
\left| \frac{|\grad u| \cdot (\log^+ \delta)^{B_2} \cdot S(\log^+ \delta )}
{\delta^{\alpha}} \right|_p, \eqno(16)
$$
where $ S(z), z \in (1,\infty) $ is positive continuous {\it slowly varying }
as $ z \to \infty $ function.\par

\vspace{3mm}

{\bf B.} Variable exponent weight estimations.\par
We consider in this section the case of bounded domain $ D. $ Let $ w_1(x), \ w_2(x) $
be two positive measurable local integrable functions (weights). \par
We denote the weight $ L_{p,w} $ average

$$
u_{w} = u_{w,D} = \frac{\int_D u(x) \ w(x) \ dx }{\int_D w(x) \ dx}
$$
and correspondence $ L_{p,w} $ norm

$$
|u|_{p,w} = \left[ \int_D |u(x)|^p \ w(x) \ dx \right]^{1/p}, \ p \ge 1.
$$
It is obtained in many publications, under some natural conditions,
 see, for example, \cite{Rjtva1}, the following
variable exponent, or the so-called $ p,q; \ p \le q $ estimations:

$$
|u - u_{w_1,D}|_{q,w_1} \le K(p,q; w_1,w_2,D) \ | \ |\grad u| \ |_{p,w_2}. \eqno(17)
$$

For example,  $ q > p(1 + \alpha) \ \Rightarrow  $

$$
\left| \frac{ u - u_{w_1,D}}{\delta^{1 + \alpha}} \right|_q \le C_{\alpha,D} \
 | q - p(1+\alpha)| ^{-1 + 1/(p(1 + \alpha)) - 1/q} \cdot
\left| \frac{ |\grad u | }{\delta^{\alpha} } \right|_{p(1+\alpha)}, \eqno(18)
$$
and the last inequality (18) is sharp. \par
 As a consequence: if
 $$
  | \grad u |/\delta^{\alpha} \in G(\psi)
 $$
 and if we denote

 $$
 \nu(q) = \nu_{\alpha,D}(q) = \inf_{p \in [1,q/(1 + \alpha)) }
 \left[ |q - p(1 + \alpha)|^{-1 + 1/(p(1 + \alpha)) - 1/q} \cdot
 \psi(p) \right],
 $$
then

$$
\left| \left| \frac{ u - u_{w_1,D}}{\delta^{1 + \alpha}} \right| \right|G(\nu)
\le C(\alpha,D) \ || \ |\grad u|/\delta^{\alpha} \ ||G(\psi). \eqno(19)
$$

\end{document}